\documentclass[12pt]{article}

\usepackage{latexsym,amsfonts,amssymb,exscale,enumerate}
\usepackage{amsmath}
\usepackage{amsthm}
 
\usepackage{graphicx}
\usepackage[all]{xy}
\SelectTips{cm}{}

\newcommand{\seq}{{\rm Seq}}
\renewcommand{\bot}{{\rm bot}}
\renewcommand{\top}{{\rm top}}
\newcommand{\Pol}{\cal{P}o\ell}
\newcommand{\ii}{ \textbf{\textit{i}}}
\newcommand{\jj}{ \textbf{\textit{j}}}
\newcommand{\jRi}{{_{\jj}R(\nu)_{\ii}}}

\newcommand{\BNC}{NH}

\newcommand{\refequal}[1]{\xy {\ar@{=}^{#1}
(-1,0)*{};(1,0)*{}};
\endxy}

\usepackage{fancyheadings}
\newcommand{\maps}{\colon}

\newcommand{\scs}{\scriptstyle}

\theoremstyle{definition}
\newtheorem{thm}{Theorem}
\newtheorem{cor}[thm]{Corollary}

\newtheorem{lem}[thm]{Lemma}
\newtheorem{rem}[thm]{Remark}
\newtheorem{prop}[thm]{Proposition}
\newtheorem{example}[thm]{Example}

\def\emph#1{{\sl #1\/}}
\let\hat=\widehat

\usepackage{bbm}
\def\N{{\mathbbm N}}
\def\Z{{\mathbbm Z}}
\def\Q{{\mathbbm Q}}

\def\nn{\notag}

\def\cal#1{\mathcal{#1}}

\def\lra{{\longrightarrow}}
\def\pmod{{\mathrm{-pmod}}}  

\def\Id{\mathrm{Id}}
\def\mc{\mathcal}
\def\mf{\mathfrak}
\def\Af{{_{\mc{A}}\mathbf{f}}}    
\newcommand{\define}{\stackrel{\mbox{\scriptsize{def}}}{=}}

\newcommand{\up}[1]{\xybox{
   (-3,-13)*{};
  (3,8)*{};
 (0,0)*{\includegraphics[scale=0.5]{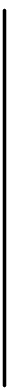}};
 (-.1,-12)*{\scs #1};
 }}

 \newcommand{\updot}[1]{\xybox{
   (-3,-13)*{};
  (3,8)*{};
 (0,0)*{\includegraphics[scale=0.5]{up.eps}};
 (-.1,-12)*{\scs #1};(0,0)*{\bullet};
 }}

\renewcommand{\sup}[1]{\xybox{
   (-3,-7)*{};
  (3,6)*{};
 (0,0)*{\includegraphics[scale=0.5]{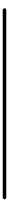}};
 (-.1,-7)*{\scs #1};
 }}

 \newcommand{\supdot}[1]{\xybox{
   (-3,-7)*{};
  (3,6)*{};
 (0,0)*{\includegraphics[scale=0.5]{short_up.eps}};
 (-.1,-7)*{\scs #1}; (0,0)*{\bullet};
 }}

\newcommand{\dcross}[2]{\xybox{
 (-6,-7)*{};
 (6,6)*{};
 (0,0)*{\includegraphics[scale=0.5]{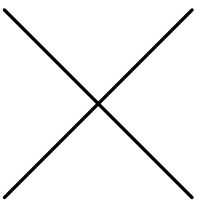}};
 (-5.1,-7)*{\scs #1};
 (4.7,-7)*{\scs #2};
}}

\newcommand{\dcrossul}[2]{\xybox{
(-2.5,2.5)*{\bullet};
 (-6,-7)*{};
 (6,6)*{};
 (0,0)*{\includegraphics[scale=0.5]{dcross.eps}};
 (-5.1,-7)*{\scs #1};
 (4.7,-7)*{\scs #2};
}}

\newcommand{\dcrossur}[2]{\xybox{
(2.5,2.5)*{\bullet};
 (-6,-7)*{};
 (6,6)*{};
 (0,0)*{\includegraphics[scale=0.5]{dcross.eps}};
 (-5.1,-7)*{\scs #1};
 (4.7,-7)*{\scs #2};
}}

\newcommand{\dcrossdl}[2]{\xybox{
(-2.5,-2.5)*{\bullet};
 (-6,-7)*{};
 (6,6)*{};
 (0,0)*{\includegraphics[scale=0.5]{dcross.eps}};
 (-5.1,-7)*{\scs #1};
 (4.7,-7)*{\scs #2};
}}

\newcommand{\dcrossdr}[2]{\xybox{
(2.5,-2.5)*{\bullet};
 (-6,-7)*{};
 (6,6)*{};
 (0,0)*{\includegraphics[scale=0.5]{dcross.eps}};
 (-5.1,-7)*{\scs #1};
 (4.7,-7)*{\scs #2};
}}

\newcommand{\twocross}[2]{\xybox{
 (-6,-13)*{};
 (6,8)*{};
 (0,0)*{\includegraphics[scale=0.5]{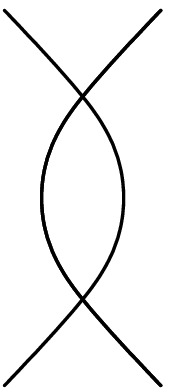}};
 (-4.1,-12)*{\scs #1};
 (3.7,-12)*{\scs #2};
}}

\newcommand{\linecrossL}[3]{\xybox{
 (-6,-13)*{};
 (6,8)*{};
 (0,0)*{\includegraphics[scale=0.5]{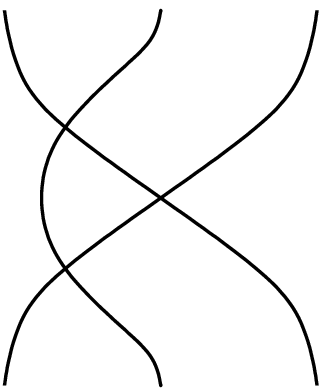}};
 (-8,-12)*{\scs #1};
 (0,-12)*{\scs #2};
 (8,-12)*{\scs #3};
}}

\newcommand{\linecrossR}[3]{\xybox{
 (-6,-13)*{};
 (6,8)*{};
 (0,0)*{\includegraphics[angle=180, scale=0.5]{line_crossL.eps}};
 (-8,-12)*{\scs #1};
 (0,-12)*{\scs #2};
 (8,-12)*{\scs #3};
}}

\title{A diagrammatic approach to categorification of quantum groups II}
      \author{ Mikhail Khovanov and Aaron D. Lauda}

\begin{document}

\date{April 15, 2008}

\maketitle

\begin{abstract}
We categorify one-half of the quantum group associated to an arbitrary Cartan datum.
\end{abstract}

\vspace{0.15in}


\noindent {\bf Cartan data and algebras $\Af$.} A \emph{Cartan datum} $(I,\cdot)$
consists of a finite set $I$ and a symmetric bilinear form on $\Z[I]$ taking
values in $\Z$, subject to conditions
\begin{itemize}
\item  $i\cdot i \in \{2, 4, 6, \dots \}$ for any $i\in I$,
\item $2\frac{i\cdot j}{i\cdot i}\in \{0, -1, -2, \dots \}$ for any $i\not= j$ in $I$.
\end{itemize}
We set $d_{ij}=-2\frac{i\cdot j}{i\cdot i}\in \N$.
To a Cartan datum assign a graph $\Gamma$ with the set of vertices $I$ and
an edge between $i$ and $j$ if and only if $i\cdot j \not= 0$.

We recall the definition of the negative half of the quantum group associated
to a Cartan datum, following~\cite{Lus4}.
Let $q_i=q^{\frac{i\cdot i}{2}}$, $[n]_i=q_i^{n-1}+q_i^{n-3}+\dots + q_i^{1-n}$,
$[n]_i!=[n]_i[n-1]_i\dots [1]_i$. Let $'\mathbf{f}$ be the free associative algebra
over $\Q(q)$ with generators $\theta_i$, $i\in I$ and denote
$\theta_i^{(n)}= \theta_i^n/[n]_i!$.  We equip $'\mathbf{f}$ with
an $\N[I]$-grading by assigning to $\theta_i$ grading $i$. The tensor square
$'\mathbf{f}\otimes {}'\mathbf{f}$ is an associative algebra with the multiplication
$$ (x_1\otimes x_2) (x_1'\otimes x_2') =q^{- |x_2|\cdot |x_1'|} x_1 x_1' \otimes x_2 x_2'$$
for homogeneous $x_1, x_2, x_1', x_2'$. There is a unique algebra homomorphism
$r: {}'\mathbf{f}\lra {}'\mathbf{f}\otimes {}'\mathbf{f}$ given on generators by
$r(\theta_i) = \theta_i \otimes 1 + 1\otimes \theta_i$.

\newpage

\begin{prop} Algebra $'\mathbf{f}$ carries a unique $\Q(v)$-bilinear form such that $(1,1)=1$
and
\begin{itemize}
\item $(\theta_i, \theta_j) = \delta_{i,j} (1-q_i^2)^{-1} $ for all $i,j\in I$,
\item $(x, y y') = (r(x), y \otimes y')$ for $x,y, y' \in {}'\mathbf{f}$,
\item $(x x', y) = (x \otimes x', r(y))$ for $x, x', y\in {}'\mathbf{f}$.
\end{itemize}
This bilinear form is symmetric.
\end{prop}

The radical $\mf{I}$ of $(,)$ is a two-sided ideal of $'\mathbf{f}$. The
bilinear form descends to a non-degenerate bilinear form on the
associative $\Q(q)$-algebra $\mathbf{f} = {}'\mathbf{f}/\mf{I}$. The
$\N[I]$-grading also descends:
$$ \mathbf{f} = \bigoplus_{\nu\in \N[I]} \mathbf{f}_{\nu}.$$
The quantum version of the Gabber-Kac theorem says the following.
\begin{prop} The ideal $\mf{I}$ is generated by the elements
$$ \sum_{a+b=d_{ij}+1} (-1)^a \theta_i^{(a)} \theta_j \theta_i^{(b)} $$
over all $i,j\in I, $ $i\not= j$.
\end{prop}
Thus, $\mathbf{f}$ is the quotient of $'\mathbf{f}$ by the so-called
quantum Serre relations
\begin{equation}\label{rels-serre}
 \sum_{a+b=d_{ij}+1} (-1)^a \theta_i^{(a)} \theta_j \theta_i^{(b)} =0.
\end{equation}

Denote by $\Af$ the $\Z[q,q^{-1}]$-subalgebra of $\mathbf{f}$ generated
by the divided powers $\theta_i^{(a)}$, over all $i\in I$ and $a\in \N$.

\vspace{0.2in}


\noindent {\bf Algebras $R(\nu)$.} As in~\cite{KL}, we consider braid-like planar
diagrams, each strand labelled by an element of $I$, and impose the following
relations
\begin{eqnarray} \label{new_eq_UUzero}
   \xy   (0,0)*{\twocross{i}{j}}; \endxy
 & = & \left\{
\begin{array}{ccl}
  0 & \qquad & \text{if $i=j$, } \\ \\
  \xy (0,0)*{\sup{i}};  (8,0)*{\sup{j}};  \endxy
  & &
 \text{if $i \cdot j=0$, }
  \\    \\
  \vcenter{ \xy
   (0,0)*{\supdot{i}};
   (8,0)*{\sup{j}};
   (-3.5,2)*{\scs d_{ij}};
  \endxy}
    \;\; + \;\;
   \vcenter{\xy  (8,0)*{\supdot{j}};(0,0)*{\sup{i}};
   (11.5,2)*{\scs d_{ji}};
   \endxy}
 & &
 \text{if $i \cdot j\not= 0$}.
\end{array}
\right.
\end{eqnarray}

\begin{eqnarray}\label{new_eq_ijslide}
  \xy  (0,0)*{\dcrossul{i}{j}};  \endxy
 \quad  =
   \xy  (0,0)*{\dcrossdr{i}{j}};   \endxy
& \quad &
   \xy  (0,0)*{\dcrossur{i}{j}};  \endxy
 \quad = \;\;
   \xy  (0,0)*{\dcrossdl{i}{j}};  \endxy \qquad \text{for $i \neq j$}
\end{eqnarray}

\begin{eqnarray}        \label{new_eq_iislide1}
 \xy  (0,0)*{\dcrossul{i}{i}}; \endxy
    \quad - \quad
 \xy (0,0)*{\dcrossdr{i}{i}}; \endxy
  & = &
 \xy (-3,0)*{\sup{i}}; (3,0)*{\sup{i}}; \endxy \\      \label{eq_iislide2}
  \xy (0,0)*{\dcrossdl{i}{i}}; \endxy
 \quad - \quad
 \xy (0,0)*{\dcrossur{i}{i}}; \endxy
  & = &
 \xy (-3,0)*{\sup{i}}; (3,0)*{\sup{i}}; \endxy
\end{eqnarray}

\begin{eqnarray}      \label{new_eq_r3_easy}
\xy  (0,0)*{\linecrossL{i}{j}{k}}; \endxy
  &=&
\xy (0,0)*{\linecrossR{i}{j}{k}}; \endxy
 \qquad \text{unless $i=k$ and $i \cdot j\not= 0$   \hspace{1in} }
\\
\xy (0,0)*{\linecrossL{i}{j}{i}}; \endxy
  &-&
\xy (0,0)*{\linecrossR{i}{j}{i}}; \endxy
 \quad = \quad
 \sum_{a=0}^{d_{ij}-1}
 \xy
 (-9,0)*{\updot{i}};
 (-6.5,3)*{\scs a};
 (0,0)*{\up{j}};
 (9,0)*{\updot{i}};
 (17,3)*{\scs d_{ij}-1-a};\endxy
 \qquad \text{if $i \cdot j\neq 0 $} \nn\\ \label{eq_r3_hard}
\end{eqnarray}

\begin{example}
For the Cartan datum $\mathbf{B}_2 = \{ i \cdot i =2, j \cdot j =4, i\cdot j =-2
\}$ we have $d_{ij}=2$, $d_{ji}=1$, and the relations involving $d_{ij}, d_{ji}$
are
\begin{equation}
  \xy   (0,0)*{\twocross{i}{j}}; \endxy
  \quad = \quad
\xy  (0,0)*{\up{i}};
   (0,5)*{\bullet};(0,-2)*{\bullet};   (8,0)*{\up{j}};  \endxy
  \;\; + \;\;
   \xy  (8,0)*{\updot{j}};
 (0,0)*{\up{i}};  \endxy \nn
\end{equation}
\begin{equation}
  \xy   (0,0)*{\twocross{j}{i}}; \endxy
  \quad = \quad
\xy  (0,0)*{\updot{j}};
   (8,0)*{\up{i}};  \endxy
  \;\; + \;\;
   \xy  (8,0)*{\up{i}}; (8,5)*{\bullet};(8,-2)*{\bullet}; (0,0)*{\up{j}};  \endxy \nn
\end{equation}
\begin{equation}
\xy (0,0)*{\linecrossL{i}{j}{i}}; \endxy
  \;\; -\;\;
\xy (0,0)*{\linecrossR{i}{j}{i}}; \endxy
 \quad = \quad
 \xy(-9,0)*{\updot{i}};(0,0)*{\up{j}};(9,0)*{\up{i}};\endxy
 \;\;+\;\;
 \xy(-9,0)*{\up{i}};(0,0)*{\up{j}};(9,0)*{\updot{i}};\endxy \nn
\end{equation}
\begin{equation}
\xy (0,0)*{\linecrossL{j}{i}{j}}; \endxy
  \;\; -\;\;
\xy (0,0)*{\linecrossR{j}{i}{j}}; \endxy
 \quad = \quad
 \xy  (-9,0)*{\up{j}}; (0,0)*{\up{i}}; (9,0)*{\up{j}};\endxy \nn
\end{equation}
\end{example}

For each $\nu\in \N[I]$ define the graded ring
  \begin{equation}
  R(\nu) \define \bigoplus_{\ii,\jj \in \seq(\nu)} \jRi ,
\end{equation}
where $\jRi$ is the abelian group of all linear combinations
of diagrams with $\bot(D)=\ii$ and $\top(D)=\jj$ modulo the relations
\eqref{new_eq_UUzero}--\eqref{eq_r3_hard} and $\seq(\nu)$ is the set of weight
$\nu$ sequences of elements of $I$. The multiplication is given
by concatenation. Degrees of the generators are
\begin{equation}
  \deg\left( \xy  (-3,0)*{\supdot{i}}; \endxy \right) = i\cdot i , \qquad \deg\left(  \xy
  (0,0)*{\dcross{i}{j}};  \endxy \right) = - i \cdot j.
\end{equation}
The rest of~\cite[Section 2.1]{KL} generalizes without difficulty to
an arbitrary Cartan datum. To define the analogue of the module $\Pol_{\nu}$ over
$R(\nu)$, we choose an orientation of each edge of $\Gamma$, then
faithfully follow the exposition in Section 2.3 of~\cite{KL},
only changing the action of $\delta_{k,\ii}$ in the last of the four cases to
$$ f \ \mapsto \  (x_k(s_k\ii)^d+x_{k+1}(s_k\ii)^{d'})(s_k f) \quad \text{if } \;\;
 i_k \longrightarrow i_{k+1}, $$
where $d=d_{i_{k+1} i_{k}}$ and $d'=d_{i_{k}i_{k+1}}$. Here notation $i_k
\longrightarrow i_{k+1}$ means that $i_k\cdot i_{k+1}\not= 0$ and this edge of
$\Gamma$ is oriented from $i_k$ to $i_{k+1}$. Proposition~2.3 in~\cite{KL} holds
for an arbitrary $(I,\cdot)$. As in~\cite[Section 2.3]{KL}, we define
${}_{\jj}B_{\ii}$, which might depend on minimal presentations of permutations in
${}_{\jj}S_{\ii}$ and gives a basis in ${}_{\jj}R(\nu)_{\ii}$. 
Corollary 2.6 in~\cite{KL}, showing that $\Pol_{\nu}$ is
a faithful graded module over $R(\nu)$, holds for an arbitrary Cartan datum and
the properties of $R(\nu)$ established in \cite[Section 2.4]{KL} generalize
without difficulty.

\vspace{0.2in}


\noindent {\bf Computations in the nilHecke ring.} In this section we slightly
enhance the graphical calculus for computations in the nilHecke ring and record
several lemmas to be used in the proof of categorified quantum Serre relations
below. We use notations from Section 2.2 of~\cite{KL}.

A box with $n$ incoming and $n$ outgoing edges and $\partial(n)$ written
inside denotes the longest divided difference $\partial_{w_0}$, the nonzero
product of $\frac{n(n-1)}{2}$ divided differences from $\{ \partial_1, \dots, \partial_{n-1}\}$:
\begin{equation}
  \xy
 (0,0)*{\includegraphics[scale=0.5]{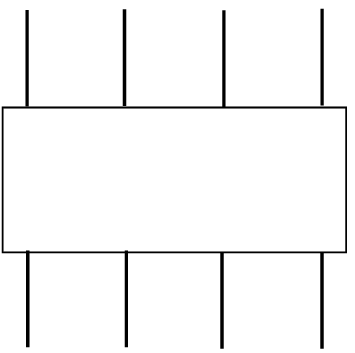}};
 (0,0)*{\partial(n)};(0,-11)*{\underbrace{\hspace{0.7in}}};
 (0,-14)*{n};
  \endxy
 \quad = \quad
   \xy
 (0,0)*{\includegraphics[scale=0.5]{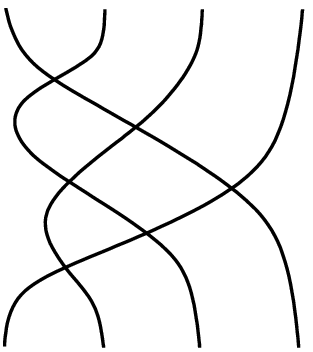}};
  \endxy \nn
\end{equation}

When this box is part of a diagram for an element of $R(\nu)$, it denotes the
corresponding element of $R(n i) \subset R(\nu)$. A box labelled $e_n$ denotes
the idempotent $ e_n = x_1^{n-1}x_2^{n-2}\dots x_{n-1}  \partial_{w_0} $.
\begin{equation}
  \xy
 (0,0)*{\includegraphics[scale=0.5]{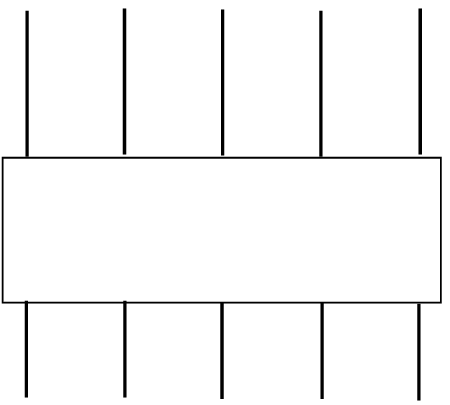}};
 (0,-1.5)*{e_{n}};
  \endxy
 \quad = \quad
 \xy
 (0,0)*{\includegraphics[scale=0.5]{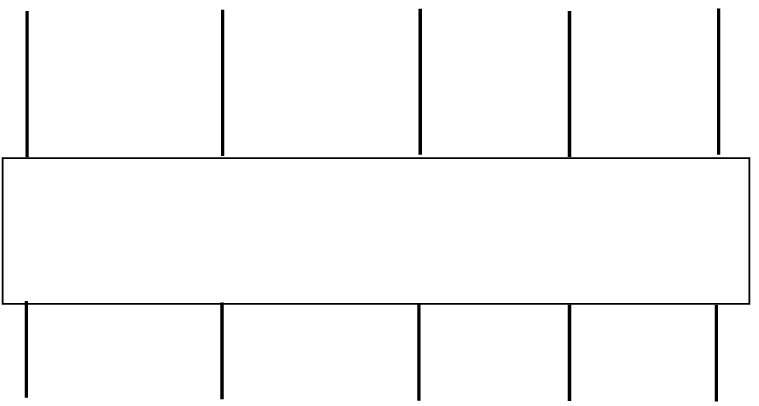}};
 (-17.7,6)*{\bullet}+(-4,2)*{\scs n-1};
 (-7.7,6)*{\bullet}+(-4,2)*{\scs n-2};
 (-2.5,5)*{\cdots};(-2.5,-7)*{\cdots};
 (2.3,6)*{\bullet}+(-2,2)*{\scs 2};
 (9.9,6)*{\bullet};(0,-1.5)*{\partial(n)};
  \endxy \nn
\end{equation}

A box labelled $e_{i,n}$ denotes the corresponding idempotent in $R(\nu)$:
\begin{equation}
  \xy
 (0,0)*{\includegraphics[scale=0.5]{c2-1.eps}};
 (0,-1.5)*{e_{i,n}};
  (-10,-12)*{\scs i};(-5,-12)*{\scs i};(0,-12)*{\cdots };
  (5,-12)*{\scs i};(10,-12)*{\scs i};
  \endxy
 \quad = \quad
 \xy
 (0,0)*{\includegraphics[scale=0.5]{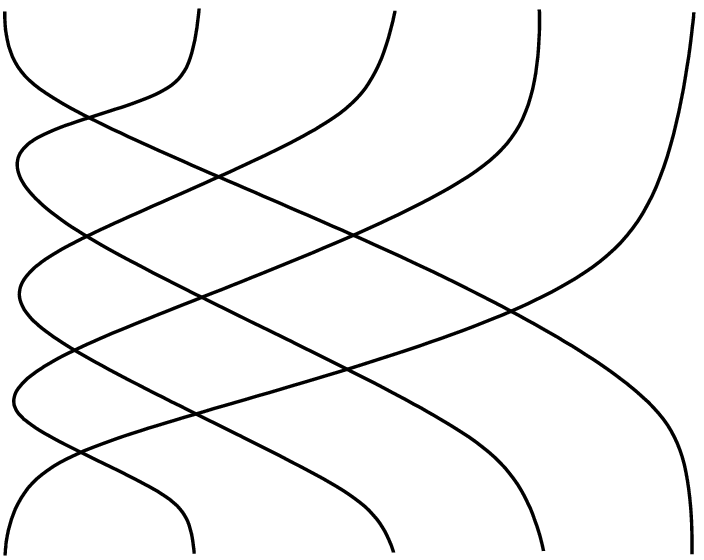}};
 (-16.7,11)*{\bullet}+(-4,1)*{\scs n-1}; (-8,11)*{\bullet}+(-4,1)*{\scs n-2};
 (1.5 ,11)*{\bullet}+(-2.2,1)*{\scs 2}; (9.5,11  )*{\bullet};
 (-18,-16)*{\scs i};(-8,-16)*{\scs i};(-2.5,-16)*{\cdots };
 (2.5,-16)*{\scs i};(10,-16)*{\scs i};(17.5,-16)*{\scs i}
  \endxy \nn
\end{equation}

\begin{rem} Similar-looking diagrams are used in the graphical calculus of
Jones-Wenzl projectors, see~\cite{KaL}, but the latter has no direct relation
to the graphical calculus in our paper.
\end{rem}

\begin{lem} We have
\begin{equation}
  \xy
 (0,0)*{\includegraphics[scale=0.5]{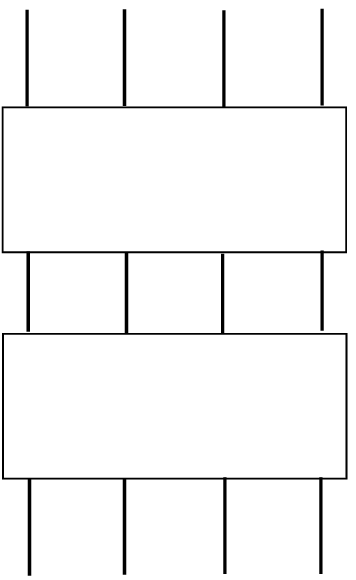}};
 (0,-5.5)*{e_{i,n}};(0,5.5)*{\partial(n)};
  \endxy
 \quad = \quad
  \xy
 (0,0)*{\includegraphics[scale=0.5]{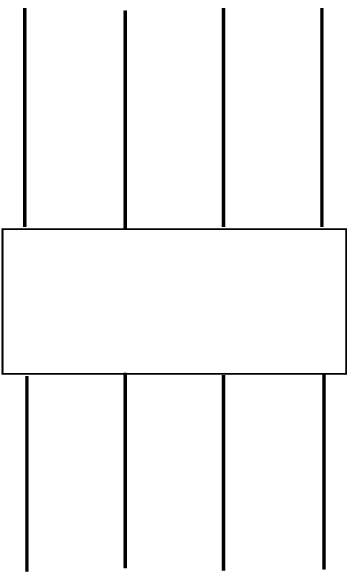}};
 (0,-0.5)*{\partial(n)};
  \endxy \nn
\end{equation}
\end{lem}

Proof is by induction on $n$:
\begin{equation}
  \xy
 (0,0)*{\includegraphics[scale=0.5]{c4-1.eps}};
 (0,-5.5)*{e_{i,n}}; (0,5.5)*{\scs \partial(n)}; \endxy
 \quad = \quad
  \xy
 (0,0)*{\includegraphics[scale=0.5]{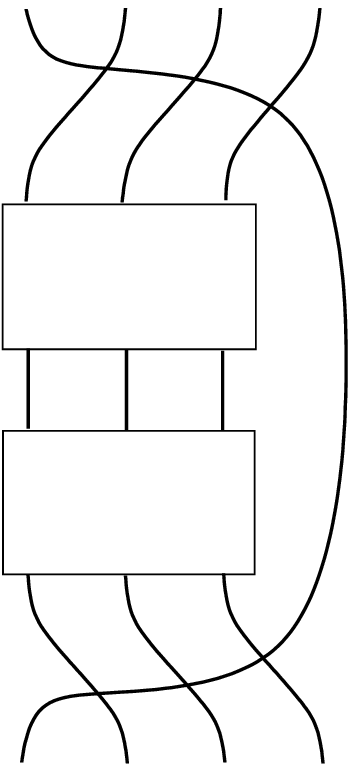}};
  (-2,-6)*{e_{i,n-1}};(-7,-12)*{\bullet};(-2,-12)*{\bullet};
  (3.2,-12)*{\bullet};(-2,5.5)*{\scs \partial(n-1)};
  \endxy
  \quad = \quad
   \xy
 (0,0)*{\includegraphics[scale=0.5]{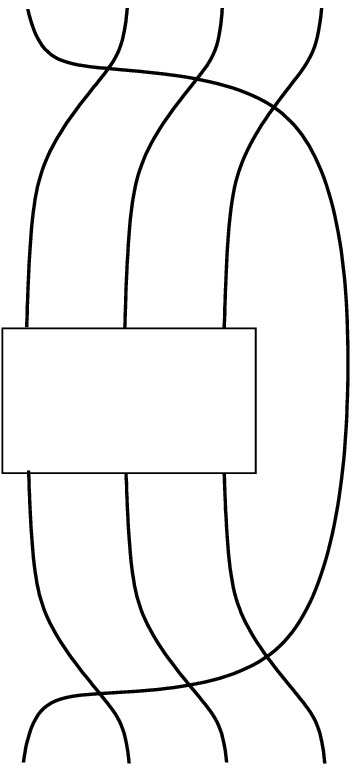}};(-2,-1)*{\scs \partial(n-1)};
 (-7,-10)*{\bullet};(-2,-10)*{\bullet}; (3,-10)*{\bullet};\endxy
  \quad = \quad
 \xy
 (1,-2)*{\includegraphics[scale=0.5]{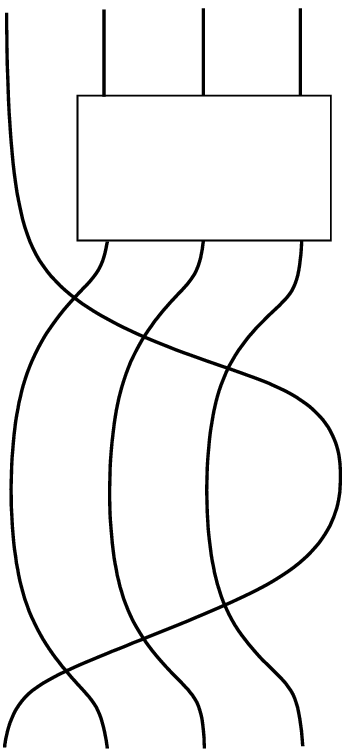}}; (3,8.5)*{\scs \partial(n-1)};
 (-7,-10)*{\bullet};(-2,-10)*{\bullet}; (3,-10)*{\bullet};\endxy \nn
\end{equation}
\begin{equation}
 = \xy
 (1,-2)*{\includegraphics[scale=0.5]{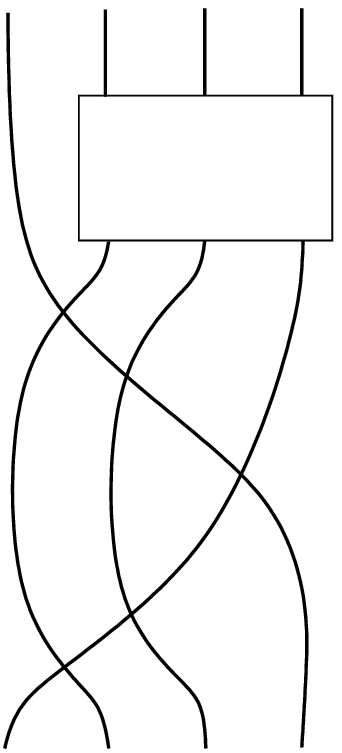}};
 (-6.8,-10)*{\bullet};(-1.8,-10)*{\bullet};(3,8.5)*{\scs \partial(n-1)}; \endxy
 \quad = \;\; \cdots \;\; = \quad \xy
 (0,0)*{\includegraphics[scale=0.5]{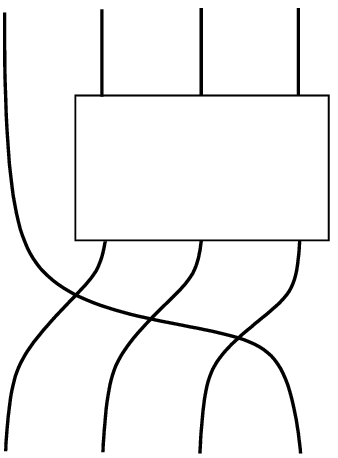}};(2.5,3)*{\scs \partial(n-1)};
  \endxy \quad = \quad
   \xy
 (0,0)*{\includegraphics[scale=0.5]{c1-1.eps}}; (0,0)*{\scs \partial(n)};\endxy \nn
\end{equation}

The first equality uses that $x_1x_2\dots x_{n-1}$ is central in the nilHecke
ring $\BNC_{n-1}$, allowing us to move these dots across $\partial(n-1)$. The
second equality is the induction hypothesis. $\square$

The lemma implies the following graphical identities:
\begin{equation} \label{eq-boxes-1}
 \xy
 (1,-2)*{\includegraphics[scale=0.5]{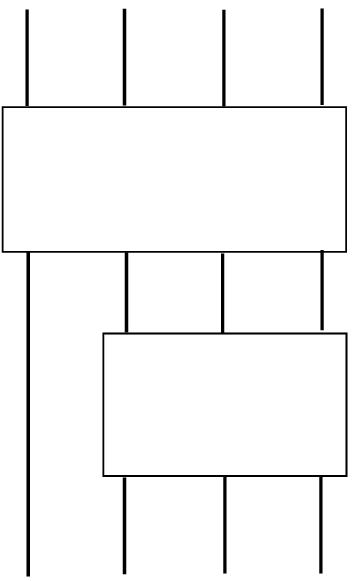}};
 (1,4)*{e_n};(4,-8)*{e_{n-1}};
\endxy
 \quad =\quad
  \xy
 (0,-2)*{\includegraphics[scale=0.5]{c4-2.eps}};
 (0,-3)*{e_n};
  \endxy
\qquad \quad
 \xy
 (1,-2)*{\includegraphics[scale=0.5]{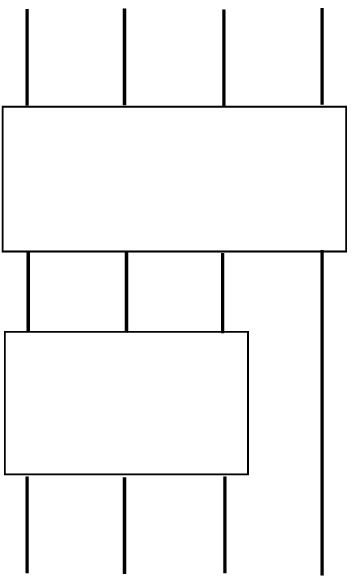}};
 (1,4)*{e_n};(-1,-8)*{e_{n-1}};
\endxy
 \quad =\quad
  \xy
 (0,-2)*{\includegraphics[scale=0.5]{c4-2.eps}};
 (0,-3)*{e_n};
  \endxy
\end{equation}
\begin{equation} \label{eq-boxes-2}
 \xy
 (-1,-2)*{\includegraphics[scale=0.5]{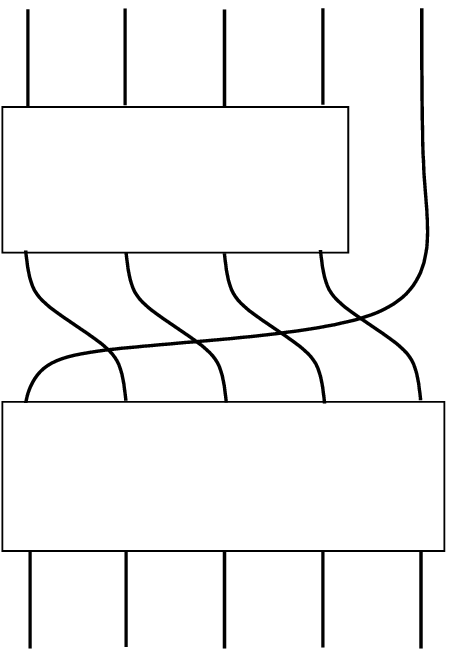}};
 (-3,5)*{e_{n-1}};(0,-10)*{e_{n}};
\endxy
 \;\; =\;\;
  \xy
 (-3,-2)*{\includegraphics[scale=0.5]{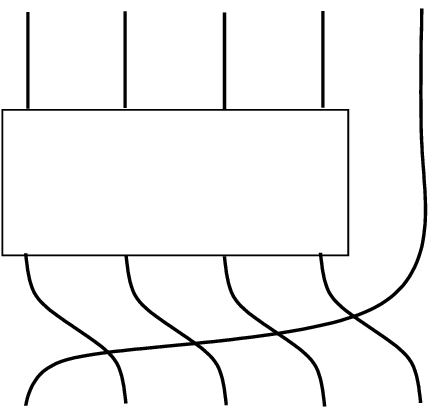}};
 (-4,-1)*{e_{n-1}};
  \endxy
\qquad \quad
 \xy
 (-1,-2)*{\reflectbox{\includegraphics[scale=0.5]{c7-1.eps}}};
 (2,5)*{e_{n-1}};(0,-10)*{e_{n}};
\endxy
 \;\; =\;\;
  \xy
 (-3,-2)*{\reflectbox{\includegraphics[scale=0.5]{c7-2.eps}}};
 (-1,-1)*{e_{n-1}};
  \endxy
\end{equation}

The following also hold:
\begin{eqnarray}
 \xy
 (1,-2)*{\includegraphics[scale=0.5]{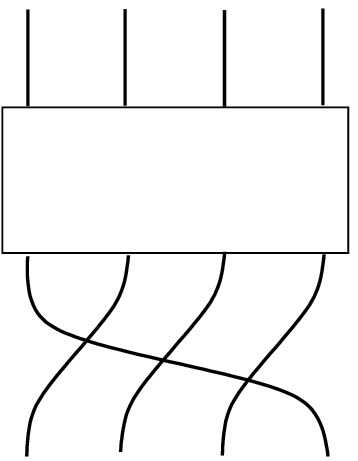}};
 (1,0)*{e_n};(-6,-6)*{\bullet}+(-2.2,0)*{\scs a};
\endxy
 & =& \left\{
\begin{array}{ccl}
  0 & &\text{if $a<n-1$,} \\ \\
    \xy
 (0,-2)*{\includegraphics[scale=0.5]{c1-1.eps}};
 (0,-2)*{e_n};
  \endxy & &\text{if $a=n-1$,}
\end{array}
 \right. \label{eq-boxes_undercross1}
\\
 \xy
 (1,-2)*{\reflectbox{ \includegraphics[scale=0.5,]{c8.eps}}};
 (1,0)*{e_n};(7.7,-6)*{\bullet}+(2.2,0)*{\scs a};
\endxy
 & =&  \left\{
\begin{array}{ccl}
  0 & &\text{if $a<n-1$,} \\ \\
   (-1)^{n-1} \;\;\xy
 (0,-2)*{\includegraphics[scale=0.5]{c1-1.eps}};
 (0,-2)*{e_n};
  \endxy & &\text{if $a=n-1$.}
\end{array}
 \right. \label{eq-boxes_undercross2}
\end{eqnarray}

\vspace{0.1in}

For each $i\in I$ the ring $R(m i)$ is isomorphic to the nilHecke ring. The
grading of a dot is now $i\cdot i$, while that of a crossing is $-i\cdot i$.
For this reason one needs to generalize the grading convention described
in~\cite[Section 2.2]{KL} and define ${}_{i,m}P$ to be the right
graded projective module   $e_{i,m}R(m i)\{-\frac{m(m-1)i\cdot i}{4}\}$,
so that the grading starts in the degree $\{-\frac{m(m-1)i\cdot i}{4}\}$.
Likewise, $P_{i,m}$ is the left graded projective module
$R(mi) \psi(e_{i,m})\{-\frac{m(m-1)i\cdot i}{4}\}$.

\vspace{0.2in}


\noindent {\bf The Grothendieck group, bilinear form and projectives.} We retain
all notations and assumptions from~\cite{KL}, working over a field $\Bbbk$,
denoting by $K_0(R(\nu))$ the Grothendieck group of the category $R(\nu)\pmod$ of
graded finitely-generated projective left $R(\nu)$-modules and forming the direct
sum
$$R = \bigoplus_{\nu\in \N[I]} R(\nu), \ \ \ K_0(R) = \bigoplus_{\nu\in \N[I]}K_0(R(\nu)).$$
Consider symmetric  $\Z[q,q^{-1}]$-bilinear form
 \begin{equation} \label{eq_bil_pair2}
  (,)\maps K_0(R(\nu)) \times K_0(R(\nu)) \longrightarrow \Z[q^{-1},q]\cdot (\nu)_q,
\end{equation}
where
\begin{equation}
  (\nu)_{q} = {\rm gdim} ({\rm Sym}(\nu)) = \prod_{i \in \Gamma}
  \left(  \prod_{a=1}^{\nu_i} \frac{1}{1-q^{a i\cdot i }} \right) ,
\end{equation}
and
\begin{equation} \label{eq-bil-ten}
([P],[Q]) = {\rm gdim}_{\Bbbk}(P^{\psi}\otimes_{R(\nu)}Q ).
\end{equation}
The character ${\rm ch}(M)$ of an $R(\nu)$-module $M$, the divided power sequences
${\rm Seqd}(\nu)$, and idempotents $1_{\ii}$ for $\ii=i_1^{(n_1)}\dots i_r^{(n_r)}
\in {\rm Seqd}(\nu)$ are defined as in~\cite[Section 2.5]{KL}. Let
$\ii != [n_1]_{i_1}! \dots [n_r]_{i_r}!$ and
$$ \langle \ii \rangle = \sum_{k=1}^r \frac{n_k(n_k-1)}{2}\cdot  \frac{i_k \cdot i_k}{2}.$$
Define graded left, respectively right, projective module
$$ P_{\ii} = R(\nu) \psi(1_{\ii}) \{ - \langle \ii \rangle \}, \ \
   {}_{\ii} P = 1_{\ii} R(\nu) \{ - \langle \ii \rangle \}.$$

\vspace{0.2in}


\noindent {\bf Quantum Serre relations.} Let
\begin{equation}
  \alpha^+_{a,b}(i,j)\quad  = \quad
  \xy
 (0,0)*{\includegraphics[scale=0.5]{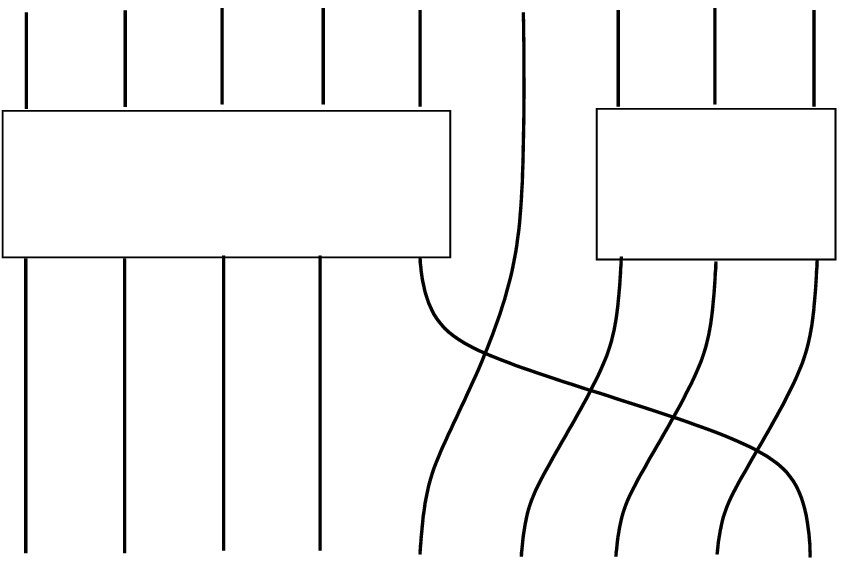}};
 (-10,5)*{e_{i,a+1}}; (15,5)*{e_{i,b-1}};
 (-20,-15.5)*{\scs i};(-15,-15.5)*{\scs i};(-10,-15.5)*{\scs i};
 (-5,-15.5)*{\scs i};(0,-15.5)*{\scs j};(5,-15.5)*{\scs i};(10,-15.5)*{\scs i};
 (15,-15.5)*{\scs i};(20,-15.5)*{\scs i};
 (-13,-19)*{\underbrace{\hspace{.7in}}};
 (12,-19)*{\underbrace{\hspace{.7in}}};
 (15,16)*{\overbrace{\hspace{.5in}}};
 (-10,16)*{\overbrace{\hspace{.85in}}};
 (-10,19)*{\scs a+1}; (15,19)*{\scs b-1};(12,-22)*{\scs b};(-13,-22)*{\scs a};
  \endxy ~,
\end{equation}
and also write $\alpha^+_{a,b}$ when $i$ and $j$ are fixed. To prove the categorified
quantum Serre relations, we assume that $a+b=d+1$, where $d=d_{ij}$.
The element $\alpha^+_{a,b}$ belongs to
\begin{equation}
  {_{i^{a+1}j i^{b-1}}R\big(j+(a+b)i\big)_{i^aj i^b}} . \nn
\end{equation}
By adding vertical lines on the left and on the right of the diagram,
$\alpha^+_{a,b}$ can be viewed, more generally, as an element of
\begin{equation}
  {_{\ii' i^{a+1}j i^{b-1}\ii''}R(\nu)_{\ii'i^aj i^b\ii''}} \nn
\end{equation}
for any sequences $\ii'$, $\ii''$ and the corresponding $\nu$. We can replace
sequences $\ii'$ and $\ii''$ by dots to simplify notation.

Left multiplication by $\alpha^+_{a,b}$ is a homomorphism of projective modules
\begin{equation}
  {_{\dots i^aji^b\dots}P} \;  \lra \; {_{\dots i^{a+1}ji^{b-1}\dots}P} .
\end{equation}
The top part of the diagram of $\alpha^+_{a,b}$ contains idempotents $e_{i,a+1}$
and $e_{i,b-1}$.  Therefore, $\alpha^+_{a,b}$ induces a homomorphism of
projective modules
\begin{equation}
  {_{\dots i^{(a)}ji^{(b)}\dots}P} \;\lra\; {_{\dots i^{(a+1)}ji^{(b-1)}\dots}P},
\end{equation}
denoted $\alpha^+_{(a,b)}$ and given by the composition
\begin{equation}
  \xymatrix@1{
  {_{\dots i^{(a)}ji^{(b)}\dots}P} \; \ar@{^{(}->}[r]
  & \; {_{\dots i^{a}ji^{b}\dots}P} \ar[r]^-{\alpha^+_{a,b}}
  & \;{_{\dots i^{(a+1)}ji^{(b-1)}\dots}P}.  }
\end{equation}
It is easy to check that $\alpha^+_{(a,b)}$ is a grading-preserving
homomorphism. Likewise, let
\begin{equation}
  \alpha^-_{a,b}(i,j) \quad  = \quad
  \xy
 (-2.5,0)*{\includegraphics[scale=0.5]{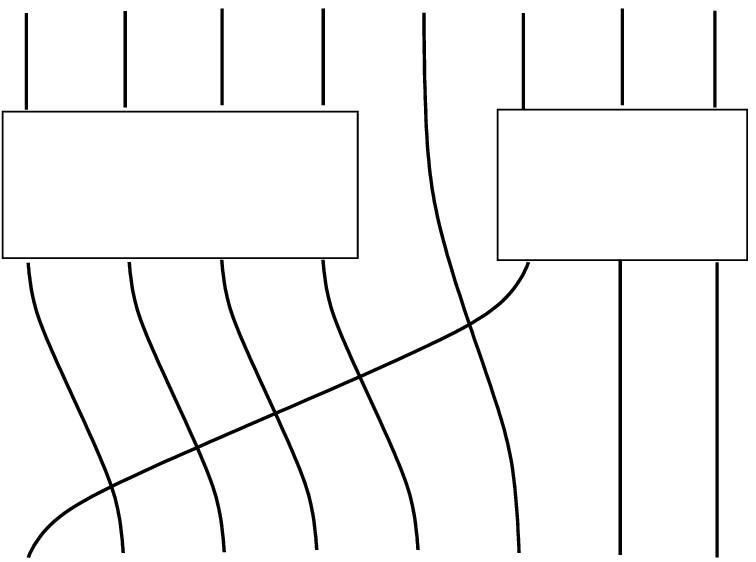}};
 (-12.5,5)*{e_{i,a-1}}; (10,5)*{e_{i,b+1}};
 (-20,-15.5)*{\scs i};(-15,-15.5)*{\scs i};(-10,-15.5)*{\scs i};
 (-5,-15.5)*{\scs i};(0,-15.5)*{\scs i};(5,-15.5)*{\scs j};(10,-15.5)*{\scs i};
 (15,-15.5)*{\scs i};
 (-10,-19)*{\underbrace{\hspace{.9in}}};
 (13,-19)*{\underbrace{\hspace{.25in}}};
 (10,16)*{\overbrace{\hspace{.5in}}};
 (-13,16)*{\overbrace{\hspace{.7in}}};
 (-13,19)*{\scs a-1}; (10,19)*{\scs b+1};(13,-22)*{\scs b};(-10,-22)*{\scs a};
  \endxy
\end{equation}
and write $\alpha^-_{a,b}$ instead of $\alpha^-_{a,b}(i,j)$ when $i$ and $j$ are
fixed.  This element of
\begin{equation}
  {_{\dots  i^{a-1}ji^{b-1}\dots}R(\nu)_{\dots i^aji^b\dots}}
\end{equation}
gives rise to a grading-preserving homomorphism of projectives
\begin{equation}
  \alpha^-_{(a,b)} \maps {_{\dots  i^{(a)}ji^{(b)}\dots}P} \;\lra\;
  {_{\dots  i^{(a-1)}ji^{(b+1)}\dots}P} .
\end{equation}

For the next few pages, denote ${_{\dots i^{(a)}ji^{(b)}\dots }P}$ by ${_{(a,b)}P}$
(recall that $a+b=d+1$, $d=d_{ij}$).
We have a diagram of projective modules and grading-preserving homomorphisms
\begin{equation} \label{eq-arrows}
  \xy
 (-35,0)*+{_{(a-1,b+1)}P}="1";(0,0)*+{_{(a,b)}P}="2";(35,0)*+{_{(a+1,b-1)}P}="3";
 (-55,0)*{}="0"; (55,0)*{}="4";
 {\ar@<1ex>^-{\alpha^+_{(a-1,b+1)}} "1";"2"};
 {\ar@<1ex>^-{\alpha^-_{(a,b)}} "2";"1"};
 {\ar@<1ex>^-{\alpha^+_{(a,b)}} "2";"3"};
 {\ar@<1ex>^-{\alpha^-_{(a+1,b-1)}} "3";"2"};
  {\ar@<1ex> "3";"4"};
 {\ar@<1ex> "4";"3"};
   {\ar@<1ex> "0";"1"};
 {\ar@<1ex> "1";"0"};
  \endxy
\end{equation}
terminating on the left at
\begin{equation}
  \xy
 (-35,0)*+{_{(0,d+1)}P}="1";
 (0,0)*+{_{(1,d)}P}="2";
 (20,0)*+{}="3";
 {\ar@<1ex>^-{\alpha^+_{(0,d+1)}} "1";"2"};
 {\ar@<1ex>^-{\alpha^-_{(1,d)}} "2";"1"};
 {\ar@<1ex> "2";"3"};
 {\ar@<1ex> "3";"2"};
  \endxy \cdots \hspace{2in}
\end{equation}
and on the right at
\begin{equation}
\cdots
  \xy
 (-20,0)*+{}="1";
 (0,0)*+{_{(d,1)}P}="2";
 (35,0)*+{_{(d+1,0)}P}="3";
 {\ar@<1ex> "1";"2"};
 {\ar@<1ex> "2";"1"};
 {\ar@<1ex>^-{\alpha^+_{(d,1)}} "2";"3"};
 {\ar@<1ex>^-{\alpha^-_{(d+1,0)}} "3";"2"};
  \endxy  .
\end{equation}

Relations~\eqref{eq-boxes-1}, \eqref{eq-boxes-2} imply
\begin{equation}
 \alpha_{a+1,b-1}^- \alpha_{a,b}^+ \quad = \quad
 \xy
 (-2,-0)*{\includegraphics[scale=0.5]{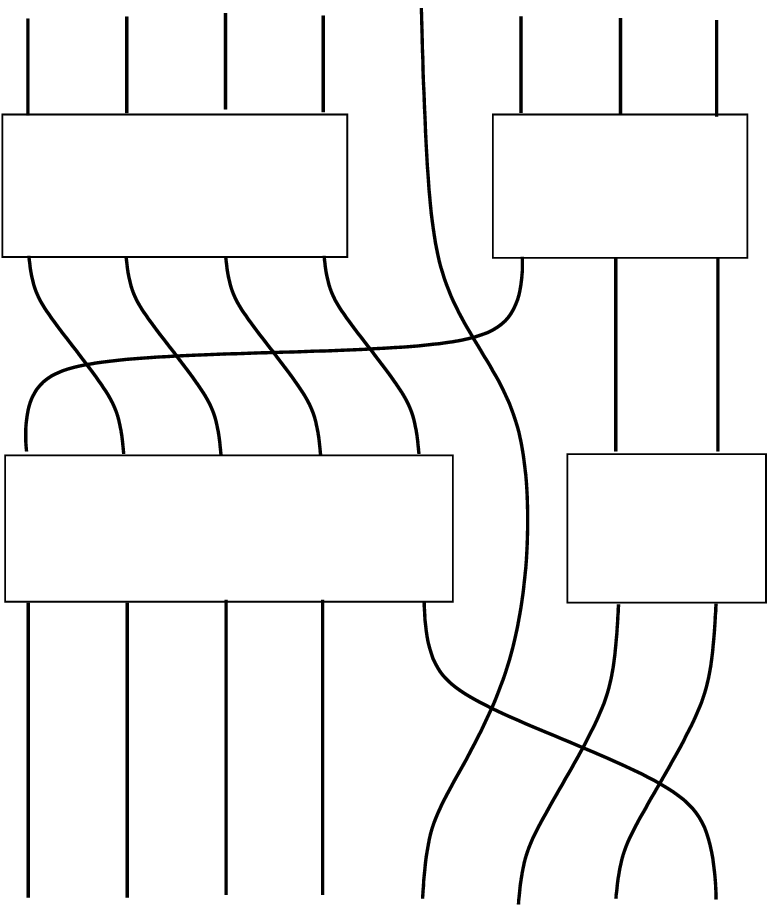}};
 (-13,14)*{e_{i,a}};(-9.6,-4)*{e_{i,a+1}};
 (10,14)*{e_{i,b}};(12,-4)*{e_{i,b-1}};
 (-20,-24.5)*{\scs i};(-15,-24.5)*{\scs i};(-10,-24.5)*{\scs i};
 (-5,-24.5)*{\scs i};(0,-24.5)*{\scs j};(5,-24.5)*{\scs i};(10,-24.5)*{\scs i};
 (15,-24.5)*{\scs i};
\endxy \nn
\end{equation}
\begin{equation}
 \quad \refequal{\eqref{eq-boxes-1}}\quad
 \xy
 (-2,0)*{\includegraphics[scale=0.5]{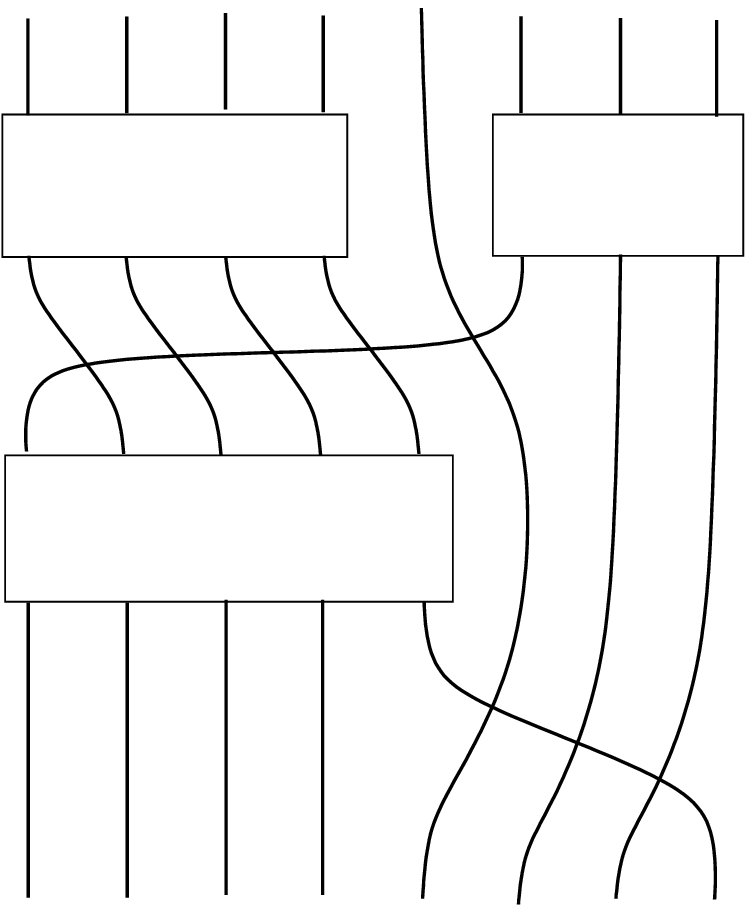}};
 (-12.5,14)*{e_{i,a}};(-9.6,-4)*{e_{i,a+1}};(10.5,14)*{e_{i,b}};
   (-20,-24.5)*{\scs i};(-15,-24.5)*{\scs i};(-10,-24.5)*{\scs i};
 (-5,-24.5)*{\scs i};(0,-24.5)*{\scs j};(5,-24.5)*{\scs i};(10,-24.5)*{\scs i};
 (15,-24.5)*{\scs i};
\endxy
 \quad \refequal{\eqref{eq-boxes-2}} \quad
 \vcenter{\xy
 (-2,-5)*{\includegraphics[scale=0.5]{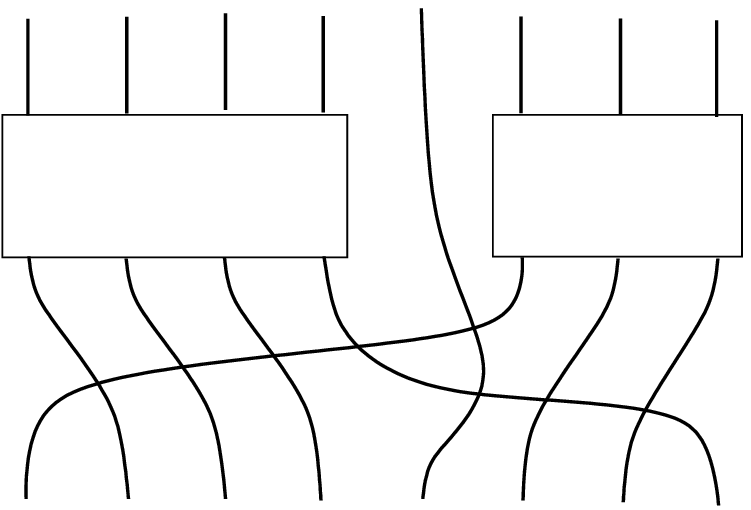}};
 (-12.5,-2)*{e_{i,a}};(10.5,-2)*{e_{i,b}};
 (0,-19.5)*{\scs j}; (-12.5,-19.5)*{\underbrace{\hspace{.65in}}};
 (10.5,-19.5)*{\underbrace{\hspace{.45in}}};(-12.5,-22)*{\scs a};
 (10.5,-22)*{\scs b};(0,12)*{};
\endxy}
\nn
\end{equation}
Furthermore,
\begin{equation}
 \alpha_{a-1,b+1}^+\alpha_{a,b}^- \quad =\quad
 \xy
 (-2,0)*{\includegraphics[scale=0.5]{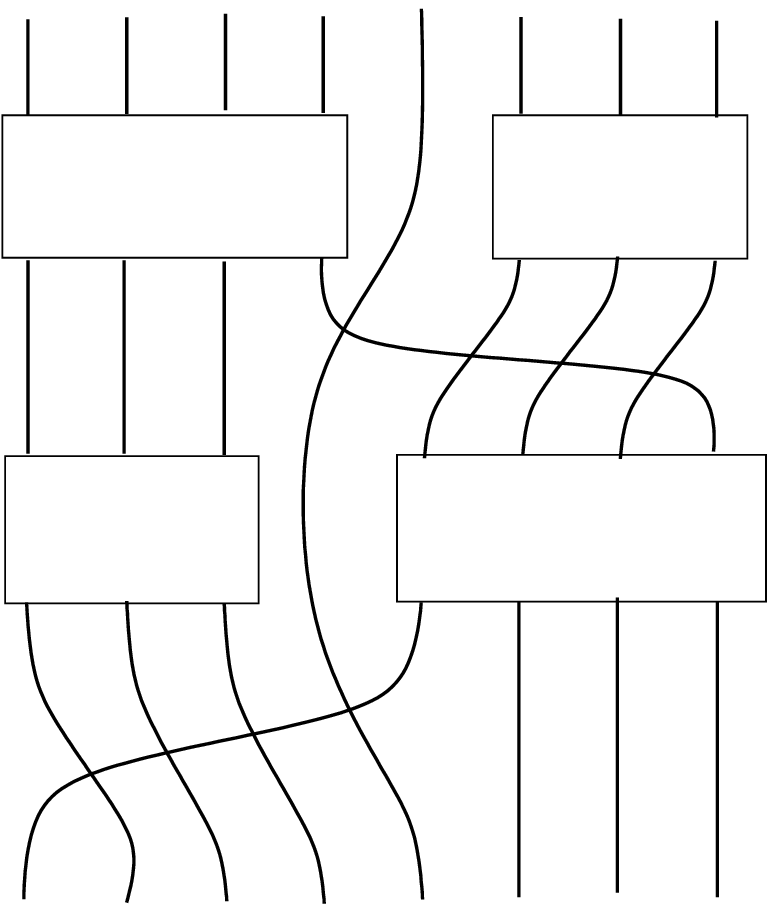}};
 (-12.5,14)*{e_{i,a}};(-14,-4)*{e_{i,a-1}};(10.5,14)*{e_{i,b}};
 (10,-4)*{e_{i,b+1}};
   (-20,-24.5)*{\scs i};(-15,-24.5)*{\scs i};(-10,-24.5)*{\scs i};
 (-5,-24.5)*{\scs i};(0,-24.5)*{\scs j};(5,-24.5)*{\scs i};(10,-24.5)*{\scs i};
 (15,-24.5)*{\scs i};
\endxy
 \quad =\quad
 \vcenter{\xy
 (-2,-10)*{\includegraphics[scale=0.5]{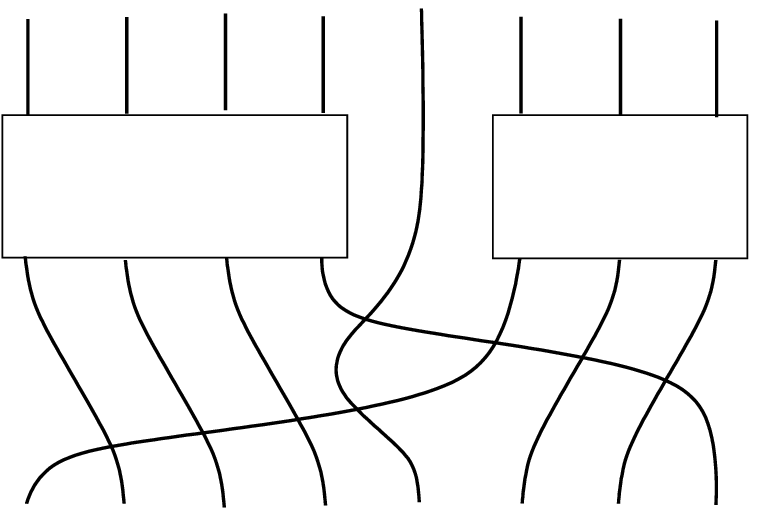}};
 (-12.5,-6)*{e_{i,a}};(10.5,-6)*{e_{i,b}};
   (-20,-24.5)*{\scs i};(-15,-24.5)*{\scs i};(-10,-24.5)*{\scs i};
 (-5,-24.5)*{\scs i};(0,-24.5)*{\scs j};(5,-24.5)*{\scs i};(10,-24.5)*{\scs i};
 (15,-24.5)*{\scs i};
\endxy}
\nn
\end{equation}

Therefore,
$$ \alpha_{a-1,b+1}^+\alpha_{a,b}^- -  \alpha_{a+1,b-1}^- \alpha_{a,b}^+ =
 (-1)^{a-1} e_{i,a} \otimes 1_j \otimes e_{i,b} ,$$
as elements of $R(\nu)$, see below:
\begin{equation}
 \alpha_{a-1,b+1}^+\alpha_{a,b}^- -  \alpha_{a+1,b-1}^- \alpha_{a,b}^+
 \;\; = \;\; \vcenter{\xy
 (-2,-10)*{\includegraphics[scale=0.5]{s2-2.eps}};
 (-12.5,-6)*{e_{i,a}};(10.5,-6)*{e_{i,b}};
   (-20,-24.5)*{\scs i};(-15,-24.5)*{\scs i};(-10,-24.5)*{\scs i};
 (-5,-24.5)*{\scs i};(0,-24.5)*{\scs j};(5,-24.5)*{\scs i};(10,-24.5)*{\scs i};
 (15,-24.5)*{\scs i};
\endxy}
\;\; - \;\;
 \vcenter{ \xy
 (-2,-10)*{\includegraphics[scale=0.5]{s1-3.eps}};
 (-12.5,-6)*{e_{i,a}};(10.5,-6)*{e_{i,b}};
   (-20,-24.5)*{\scs i};(-15,-24.5)*{\scs i};(-10,-24.5)*{\scs i};
 (-5,-24.5)*{\scs i};(0,-24.5)*{\scs j};(5,-24.5)*{\scs i};(10,-24.5)*{\scs i};
 (15,-24.5)*{\scs i};
\endxy}
\nn
\end{equation}
\begin{eqnarray}
\quad &\refequal{\eqref{eq_r3_hard}} & \sum_{c=0}^{d-1}\;\;
 \vcenter{\xy
 (-2,-10)*{\includegraphics[scale=0.5]{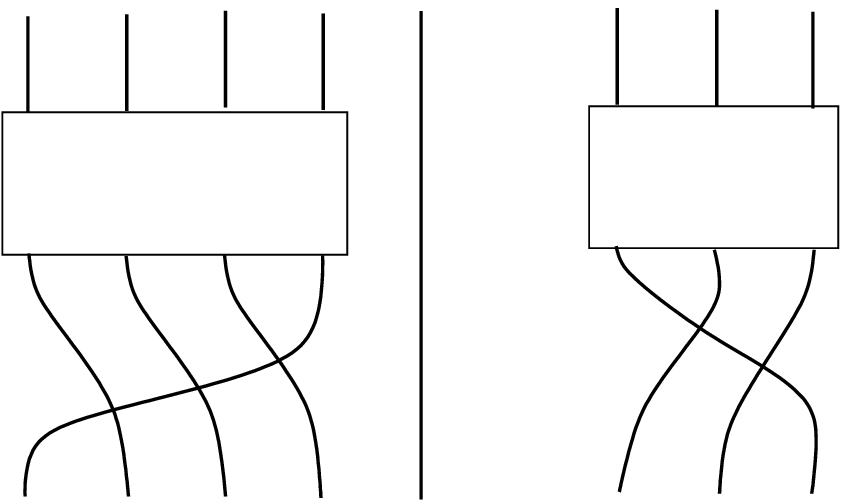}};
 (-14.5,-6.5)*{e_{i,a}};(13,-6.5)*{e_{i,b}};
   (-22.5,-24.5)*{\scs i};(-17.5,-24.5)*{\scs i};(-12.5,-24.5)*{\scs i};
 (-7.5,-24.5)*{\scs i};(-2.4,-24.5)*{\scs j};(7.5,-24.5)*{\scs i};(12.5,-24.5)*{\scs i};
 (17.5,-24.5)*{\scs i};
 (-7,-13)*{\bullet}+(2.2,0)*{\scs c};
 (10.5,-13)*{\bullet}+(-6,0)*{\scs d-1-c};
\endxy}
 \quad \refequal{\eqref{eq-boxes_undercross1},\eqref{eq-boxes_undercross2}} \quad
  \vcenter{\xy
 (-2,-10)*{\includegraphics[scale=0.5]{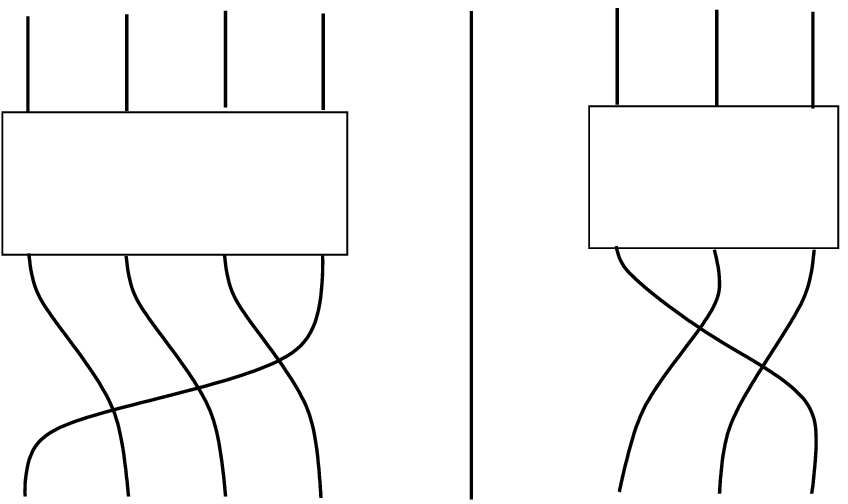}};
 (-14.5,-6.5)*{e_{i,a}};(13,-6.5)*{e_{i,b}};
   (-22.5,-24.5)*{\scs i};(-17.5,-24.5)*{\scs i};(-12.5,-24.5)*{\scs i};
 (-7.5,-24.5)*{\scs i};(-0,-24.5)*{\scs j};(7.5,-24.5)*{\scs i};(12.5,-24.5)*{\scs i};
 (17.5,-24.5)*{\scs i};
 (-7,-13)*{\bullet}+(4,0)*{\scs a-1};
 (10.5,-13)*{\bullet}+(-4,0)*{\scs b-1};
\endxy} \nn \\
\quad &\refequal{\eqref{eq-boxes_undercross1},\eqref{eq-boxes_undercross2}} &
(-1)^{a-1}\; \vcenter{\xy
 (-2.5,-13)*{\includegraphics[scale=0.5]{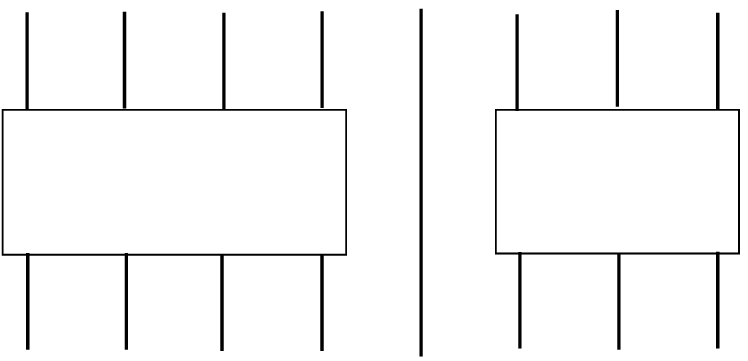}};
 (-12.5,-13)*{e_{i,a}};(11,-13.5)*{e_{i,b}};
    (-20,-24.5)*{\scs i};(-15,-24.5)*{\scs i};(-10,-24.5)*{\scs i};
 (-5,-24.5)*{\scs i};(0,-24.5)*{\scs j};(5,-24.5)*{\scs i};(10,-24.5)*{\scs i};
 (15,-24.5)*{\scs i};
\endxy}  ~,
 \nn
\end{eqnarray}
where $d=a+b-1$. Consequently,
\[\alpha_{(a-1,b+1)}^+\alpha_{(a,b)}^- -  \alpha_{(a+1,b_1)}^- \alpha_{(a,b)}^+ =
  (-1)^{a-1}\cdot \Id , \]
as endomorphisms of the projective module ${_{(a,b)}P}$, since $ e_{i,a} \otimes
1_j \otimes e_{i,b}$ acts by the identity on ${}_{(a,b)}P$.

Likewise,
\begin{eqnarray}
 \alpha^-_{1,d}\alpha^+_{0,d+1} &= &
 \vcenter{\xy
 (0.5,0)*{\includegraphics[scale=0.5]{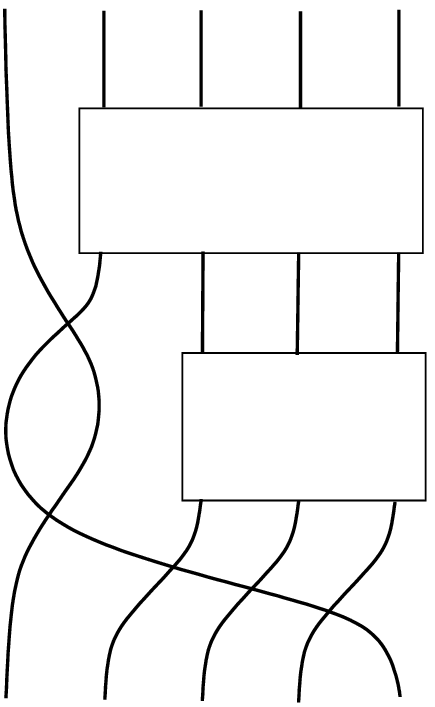}};
 (3,9)*{e_{i,d+1}};(5,-4)*{e_{i,d}};
 (-10.5,-19.5)*{\scs j};
 (-5,-19.5)*{\scs i};(0,-19.5)*{\scs i};(5,-19.5)*{\scs i};(10,-19.5)*{\scs i};
\endxy}
\quad \refequal{\eqref{new_eq_ijslide}} \quad
 \vcenter{\xy
 (0.5,0)*{\includegraphics[scale=0.5]{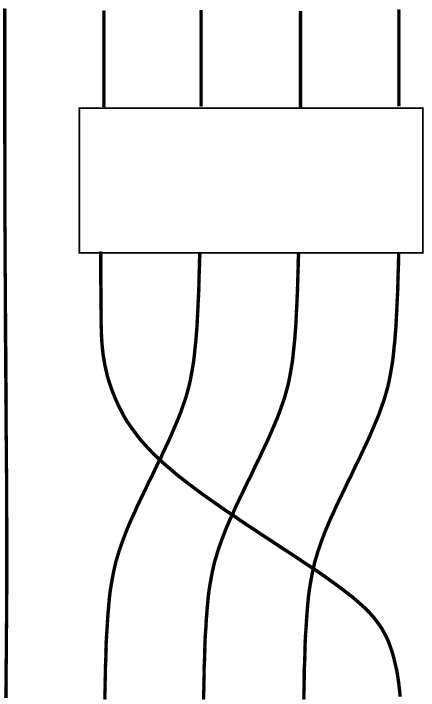}};
 (3,9)*{e_{i,d+1}};
 (-10.5,-19.5)*{\scs j};
 (-5,-19.5)*{\scs i};(0,-19.5)*{\scs i};(5,-19.5)*{\scs i};(10,-19.5)*{\scs i};
 (-10,-1)*{\bullet}+(-2.5,1)*{\scs d'};
\endxy}
\;\;+\;\; \vcenter{\xy
 (0.5,0)*{\includegraphics[scale=0.5]{s4-2.eps}};
 (3,9)*{e_{i,d+1}};
 (-10.5,-19.5)*{\scs j};
 (-5,-19.5)*{\scs i};(0,-19.5)*{\scs i};(5,-19.5)*{\scs i};(10,-19.5)*{\scs i};
 (-5,-1)*{\bullet}+(-2.5,0)*{\scs d};
\endxy} \nn \\
\quad &\refequal{\eqref{eq-boxes_undercross1}}& \vcenter{\xy
 (0.5,0)*{\includegraphics[scale=0.5]{s4-2.eps}};
 (3,9)*{e_{i,d+1}};
 (-10.5,-19.5)*{\scs j};
 (-5,-19.5)*{\scs i};(0,-19.5)*{\scs i};(5,-19.5)*{\scs i};(10,-19.5)*{\scs i};
 (-5,-1)*{\bullet}+(-2.5,0)*{\scs d};
\endxy}
\quad \refequal{\eqref{eq-boxes_undercross1}}  \quad \vcenter{\xy
 (0.5,-7)*{\includegraphics[scale=0.5]{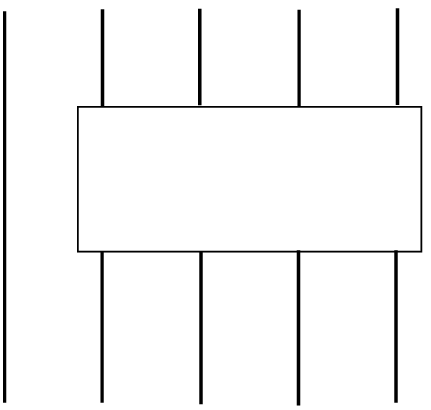}};
 (3,-6)*{e_{i,d+1}};
 (-10.5,-19.5)*{\scs j};
 (-5,-19.5)*{\scs i};(0,-19.5)*{\scs i};(5,-19.5)*{\scs i};(10,-19.5)*{\scs i};
\endxy}
\quad =\quad 1_j \otimes e_{i,d+1} \  , \nn
\end{eqnarray}
where $d'=d_{ji}$, and
$$ \alpha^-_{(1,d)}\alpha^+_{(0,d+1)} = \Id,$$
as endomorphisms of ${}_{(0,d+1)}P$. A similar computation shows that
$$ \alpha^+_{d,1}\alpha^-_{d+1,0} = (-1)^d  e_{i,d+1}\otimes 1_j, $$
as elements of $R(\nu)$, and
$$ \alpha^+_{(d,1)}\alpha^-_{(d+1,0)} = (-1)^d\cdot   \Id , $$
as endomorphisms of ${}_{(d+1,0)}P$.

\begin{prop}\label{serre-right}
 For each $i, j\in I, i\not= j$ there are isomorphisms of graded
right projective modules
$$ \bigoplus_{a=0}^{\lfloor \frac{d+1}{2} \rfloor}
 {}_{\dots i^{(2a)}j i^{(d+1-2a)} \dots}P \cong
   \bigoplus_{a=0}^{\lfloor \frac{d}{2} \rfloor}
 {}_{\dots i^{(2a+1)}j i^{(d-2a)} \dots}P  .
$$
\end{prop}
\begin{proof} When $i\cdot j=0$ the isomorphism reads
$$ {}_{\dots ji \dots} P \cong {}_{\dots ij \dots} P$$
and is given by left multiplication by the $ij$ intersection. When $i\cdot j <0$,
earlier computations show that the maps
\[ \bigoplus_{a=0}^{\lfloor
\frac{d+1}{2} \rfloor} {}_{(2a, d+1-2a)} P \xymatrix@1{
 \ar@<.8ex>[rr]^-{\alpha'} &&
  \ar@<.8ex>[ll]^-{\alpha''}
  } \bigoplus_{a=0}^{\lfloor \frac{d}{2} \rfloor}  {}_{(2a+1, d-2a)}P
\]
given by
\begin{eqnarray*}
\alpha' & = & \sum_{a=0}^{\lfloor \frac{d}{2} \rfloor}\alpha^+_{(2a,d+1-2a)} +
  \sum_{a=0}^{\lfloor \frac{d+1}{2} \rfloor}\alpha^-_{(2a,d+1-2a)},  \\
\alpha'' & = & \sum_{a=0}^{\lfloor \frac{d+1}{2} \rfloor}\alpha^-_{(2a+1,d-2a)} -
 \sum_{a=0}^{\lfloor \frac{d}{2} \rfloor}\alpha^+_{2a+1,d-2a} ,
\end{eqnarray*}
are mutually-inverse isomorphisms, implying the proposition. Maps $\alpha', \alpha''$
together are given by summing over all arrows in the diagram~\eqref{eq-arrows},
with every fourth arrow appearing with the minus sign.
\end{proof}

\begin{cor}\label{serre-left} For each $i, j\in I, i\not= j$ there are isomorphisms of graded
left projective modules
$$ \bigoplus_{a=0}^{\lfloor \frac{d+1}{2} \rfloor}
 P_{\dots i^{(2a)}j i^{(d+1-2a)} \dots} \cong
   \bigoplus_{a=0}^{\lfloor \frac{d}{2} \rfloor}
 P_{\dots i^{(2a+1)}j i^{(d-2a)} \dots}
$$
\end{cor}

Proposition~\ref{serre-right} and Corollary~\ref{serre-left} generalize
Proposition 2.13 in~\cite{KL} and can be considered a categorification of the
quantum Serre relations. Corollaries 2.14 and 2.15 of~\cite{KL}, establishing
quantum Serre relations for the characters of any $M\in R(\nu)\mathrm{-mod}$,
generalize to an arbitrary Cartan datum in the same way.

\vspace{0.2in}


\noindent {\bf Grothendieck group as the quantum group.} Induction and
restriction functors for inclusions $R(\nu)\otimes R(\nu')\subset R(\nu+\nu')$
turn $K_0(R)$ into a twisted bialgebra, and all results of~\cite[Section 2.6]{KL}
remain valid for an arbitrary Cartan datum. As in~\cite[Section 3.1]{KL} we
define a homomorphism of twisted bialgebras
$$ \gamma \ : \  \Af \lra K_0(R)$$
which takes the product of divided powers $\theta_{\ii} =
\theta_{i_1}^{(n_1)}\dots \theta_{i_r}^{(n_r)}$ to $[P_{\ii}]$, where $\ii =
i_{i_1}^{(n_1)}\dots i_{i_r}^{(n_r)}$. Homomorphism $\gamma$ intertwines the
bilinear forms on $\Af$ and $K_0(R)$,
$$ (x, y) = (\gamma(x), \gamma(y)), \ \ \  x,y\in \Af. $$
Due to the quantum Gabber-Kac theorem, this homomorphism is injective.
Surjectivity of $\gamma$ follows from the arguments identical to those given
in~\cite[Section 3.2]{KL}, which, in turn, were adopted from~\cite[Section 5]{KleBook}.
Alternatively, the arguments could be adopted from~\cite{GV}
and~\cite{Vaz1}; we settled on using a single source. We obtain

\begin{thm} $\gamma : \Af \lra K_0(R)$ is an isomorphism of $\N[I]$-graded twisted
bialgebras.
\end{thm}
This theorem holds without any restrictions on the Cartan datum and on the ground
field $\Bbbk$ over which $R(\nu)$ is defined. All other results and observations
of Sections 3.2 and 3.3 of~\cite{KL} extend to the general case as well.  The
cyclotomic quotients of $R(\nu)$ described in~\cite[Section 3.4]{KL} generalize
to an arbitrary Cartan datum.

It would be interesting to relate
our construction to Lusztig's geometric realization of $U^-$ in the non-simply laced
case~\cite{Lus4} and to Brundan-Kleshchev's categorification~\cite{BK},
\cite{KleBook} of $U^-_{q=1}$  in the affine Dynkin case $A^{(2)}_n$.

\vspace{0.2in}


\noindent {\bf A multi-grading.} For every pair $(i,j)$ of vertices of $\Gamma$,
algebras $R(\nu)$ can be equipped with an additional grading, by assigning degrees
$-1$ and $1$ to the $ij$ and $ji$ crossings, respectively,
\[
  \deg\left(\;\xy  (0,0)*{\dcross{i}{j}};  \endxy \;\right) =-1\ ,  \qquad
  \deg\left(\;\xy  (0,0)*{\dcross{j}{i}};  \endxy \;\right) =1 \ ,
\]
and degree $0$ to all other diagrammatic generators of $R(\nu)$. These gradings
are independent, and together with the principal grading, introduced above, make
$R(\nu)$ into a multi-graded ring (with $\frac{n(n-1)}{2}+1$ independent gradings
where $n=|{\rm Supp}(\nu)|$). The direct sum of the categories of multi-graded
finitely-generated projective left $R(\nu)$-modules, over all $\nu\in \N[I]$,
categorifies a multi-parameter deformation~\cite{OY}, \cite{Resh} of the
quantum universal enveloping algebra $U^-$, the quotient of the free associative
algebra on $\theta_i$, $i\in I$, by the relations
\begin{equation}\label{rels-multiserre}
 \sum_{a+b=d_{ij}+1} (-1)^a q_{ij}^a\theta_i^{(a)} \theta_j \theta_i^{(b)} =0,
\end{equation}
where $q_{ij}$ are formal variables subject to conditions $q_{ij}q_{ji}=1$.

\vspace{0.2in}


\noindent {\bf Modifications in the simply-laced case.} This section explains how
to deform algebras $R(\nu)$ in the simply-laced case  so that the main
results of~\cite{KL} will hold for the modified algebras.
These deformations can be nontrivial only when the
graph has cycles. As in~\cite{KL}, we start with an unoriented graph $\Gamma$
without loops and multiple edges. Next, fix an orientation of each edge of
$\Gamma$, work over a base field $\Bbbk$, and, for each oriented edge
$i\longrightarrow j$, choose two invertible elements $\tau_{ij}$ and $\tau_{ji}$
in $\Bbbk$. Denote such a datum $\{{\rm orientations}, \;{\rm invertible}\; {\rm
elements}\}$ by $\tau$.

For each $\nu\in\N[I]$ consider $\Bbbk$-vector space $\Pol_{\nu}$
defined as in~\cite{KL}. This space is the sum of polynomial rings
in $|\nu|$ variables, over all sequences in $\seq(\nu)$. Define $R_{\tau}(\nu)$
to be the endomorphism algebra of $\Pol_{\nu}$ generated by
the endomorphisms $1_{\ii},$  $x_{k,\ii}$, $\delta_{k,\ii}$, over all
possible $k$ and $\ii$, with the action as in~\cite[Section 2.3]{KL}, with
the only difference being the action of $\delta_{k,\ii}$ in the last of the
four cases:
$$ f \ \mapsto \  (\tau_{i_k i_{k+1}} x_{k+1}(s_k\ii)-
    \tau_{i_{k+1}i_k}x_k(s_k\ii))(s_k f) \quad \text{if } \;\;
     i_k \longrightarrow i_{k+1}, $$
    instead of
$$ f \ \mapsto \  (x_k(s_k\ii)+ x_{k+1}(s_k\ii))(s_k f) \quad
\text{if } \;\;  i_k \longrightarrow i_{k+1}. $$
The algebra $R_{\tau}(\nu)$ has a diagrammatic description
similar to that of $R(\nu)$, with the following defining relations
\begin{eqnarray} \label{eq_UUzero.2}
   \xy   (0,0)*{\twocross{i}{j}}; \endxy
 & = & \left\{
\begin{array}{ccc}
  0 & \qquad & \text{if $i=j$, } \\ \\
  \xy (0,0)*{\sup{i}};  (8,0)*{\sup{j}};  \endxy
  & &
 \text{if $i \cdot j=0$, }
  \\    \\
 \tau_{ij} \xy (0,0)*{\supdot{i}};   (8,0)*{\sup{j}};  \endxy
  \quad -\quad \tau_{ji}
   \xy  (8,0)*{\supdot{j}};  (0,0)*{\sup{i}};  \endxy
 & &
 \text{if $i \longrightarrow j$, }
 \\    \\
  \tau_{ji}\xy  (8,0)*{\supdot{j}};  (0,0)*{\sup{i}};  \endxy
   \quad -\quad \tau_{ij}
  \xy  (0,0)*{\supdot{i}};   (8,0)*{\sup{j}};  \endxy
   & &
 \text{if $i \longleftarrow j$, }
\end{array}
\right.
\end{eqnarray}

\begin{eqnarray}\label{eq_ijslide.2}
  \xy  (0,0)*{\dcrossul{i}{j}};  \endxy
 \quad  = \;\;
   \xy  (0,0)*{\dcrossdr{i}{j}};   \endxy
& \quad &
   \xy  (0,0)*{\dcrossur{i}{j}};  \endxy
 \quad = \;\;
   \xy  (0,0)*{\dcrossdl{i}{j}};  \endxy \qquad \text{for $i \neq j$}
\end{eqnarray}

\begin{eqnarray}        \label{eq_iislide1.2}
 \xy  (0,0)*{\dcrossul{i}{i}}; \endxy
    \quad - \quad
 \xy (0,0)*{\dcrossdr{i}{i}}; \endxy
  & = &
 \xy (-3,0)*{\sup{i}}; (3,0)*{\sup{i}}; \endxy \\      \label{eq_iislide2.2}
  \xy (0,0)*{\dcrossdl{i}{i}}; \endxy
 \quad - \quad
 \xy (0,0)*{\dcrossur{i}{i}}; \endxy
  & = &
 \xy (-3,0)*{\sup{i}}; (3,0)*{\sup{i}}; \endxy
\end{eqnarray}

\begin{eqnarray}      \label{eq_r3_easy.2}
\xy  (0,0)*{\linecrossL{i}{j}{k}}; \endxy
  &=&
\xy (0,0)*{\linecrossR{i}{j}{k}}; \endxy
 \qquad \text{unless $i=k$ and $i \cdot j=-1$   \hspace{1in} }
\\                   \label{eq_r3_hard.2}
\xy (0,0)*{\linecrossL{i}{j}{i}}; \endxy
  &-&
\xy (0,0)*{\linecrossR{i}{j}{i}}; \endxy
 \quad = \ \ \  \tau_{ij}\quad
 \xy  (-9,0)*{\up{i}};  (0,0)*{\up{j}}; (9,0)*{\up{i}}; \endxy
 \qquad \text{if $i \longrightarrow j$ }
\\
  \label{eq_r3_ashard.2}
\xy (0,0)*{\linecrossL{i}{j}{i}}; \endxy
  &-&
\xy (0,0)*{\linecrossR{i}{j}{i}}; \endxy
 \quad = \ - \tau_{ij}\quad
 \xy  (-9,0)*{\up{i}};  (0,0)*{\up{j}}; (9,0)*{\up{i}}; \endxy
 \qquad \text{if $i \longleftarrow j$ }
\end{eqnarray}

Reverse the orientation of a single edge $i-j$ and change $\tau_{ij}$ to $-\tau_{ij}$
and $\tau_{ji}$ to $-\tau_{ji}$. Denote the new datum by $\tau'.$ Algebras
$R_{\tau}(\nu)$ and $R_{\tau'}(\nu)$ are isomorphic via a map which is the
identity on diagrams. This way, the study of $R_{\tau}(\nu)$ reduces to the case of
any preferred orientation of $\Gamma$. Rescaling one of the two possible types
of the $ij$ crossing by $\lambda\in \Bbbk$ changes $\tau_{ij}$ to $\lambda\tau_{ij}$ and
$\tau_{ji}$ to $\lambda\tau_{ji}$ while keeping the rest of the data fixed. We see that
$R_{\tau}(\nu)$ depends only on products $\tau_{ij}\tau_{ji}^{-1}$, over all edges
of $\Gamma$, via non-canonical isomorphisms.  Rescalings of $ii$ crossings and dots
further reduce the number of parameters to the rank of the first homology group of $\Gamma$.
When graph $\Gamma$ is a forest (has no cycles), algebras $R_{\tau}(\nu)$
are all isomorphic to $R(\nu)$ via rescaling of generators. When
$\Gamma$ has a single cycle, rescaling of generators reduces this family
of algebras to a one-parameter family, with the parameter taking values in
$\Bbbk^{\ast}$. It is likely that $R_{\tau}(\nu)$ has a description via
equivariant convolution algebras in Lusztig's geometrization~\cite{Lus4} of $U^-$
when all $\tau_{ij}=1$ (compare with Conjecture~1.2 in~\cite{KL}).

Form
$$ R_{\tau} = \bigoplus_{\nu\in \N[I]} R_{\tau}(\nu).$$
The Grothendieck group $K_0(R_{\tau})$ of the category of finitely-generated
graded left projective modules can be naturally identified with the integral
version $\Af$ of $U^-$. All other essential constructions and results
of~\cite{KL} generalize from
$R(\nu)$ to algebras $R_{\tau}(\nu)$ in a straightforward fashion.

\vspace{0.2in}

\noindent {\bf Modifications in the general case.} Rings $R(\nu)$ associated to
an arbitrary Cartan datum admit similar modifications that depend on choosing an
orientation of $\Gamma$ and invertible elements $\tau_{ij},$  $\tau_{ji}$ of the
ground field $\Bbbk$ for each oriented edge $i\longrightarrow j$. The key point
is the change in the definition of the endomorphism algebra, making
$\delta_{k,\ii}$ act by
$$ f \ \mapsto \  (\tau_{i_k i_{k+1}} x_{k+1}(s_k\ii)^{d'}-
    \tau_{i_{k+1}i_k}x_k(s_k\ii)^d)(s_k f) \quad \text{if } \;\;
     i_k \longrightarrow i_{k+1} $$
in the last of the four cases, with $d=d_{i_{k+1}i_k}$ and $d'=d_{i_ki_{k+1}}$.
Our proof of categorified quantum Serre relations for $R(\nu)$ requires only minor
changes in the general case of $R_{\tau}(\nu)$. Everything else generalizes as well.

\vspace{0.2in}

\noindent {\bf Acknowledgments.} M.K. was fully supported by the IAS and the NSF
grants DMS--0635607 and  DMS-0706924  while working on this paper.



\vspace{0.1in}

\noindent
M.K.: {\sl \small School of Mathematics, Institute for Advanced Study, Princeton, NJ 08540, and}
{ \sl \small Department of Mathematics, Columbia University, New York, NY 10027}
\noindent
  {\tt \small email: khovanov@math.columbia.edu}

\vspace{0.1in}

\noindent
A.L:  { \sl \small Department of Mathematics, Columbia University, New York, NY 10027}
\noindent
  {\tt \small email: lauda@math.columbia.edu}

\end{document}